\documentclass[a4paper,10pt]{amsart}
\usepackage{amssymb,amsthm,amsmath,amscd}
\usepackage[dvips]{graphicx}
\usepackage{times,mathptmx} 
\newcommand{\Z}{\mathbf{Z}}
\newcommand{\R}{\mathbf{R}}
\newcommand{\C}{\mathbf{C}}
\newcommand{\HS}{\mathtt{S}}
\newcommand{\SO}{\operatorname{SO}}
\newcommand{\lk}{\operatorname{lk}}
\newcommand{\Int}{\operatorname{Int}}
\newcommand{\Emb}{\operatorname{Emb}}
\newcommand{\co}{\colon\thinspace}
\newcommand{\conn}{\thinspace\sharp\thinspace}
\newcommand{\bconn}{\thinspace\natural\thinspace}

\newtheorem{theorem}{Theorem}[section]
\newtheorem{corollary}[theorem]{Corollary}
\newtheorem{lemma}[theorem]{Lemma}

\newtheorem{proposition}[theorem]{Proposition}
\theoremstyle{definition}
\newtheorem{definition}[theorem]{Definition}
\newtheorem{remark}[theorem]{Remark}
\numberwithin{equation}{section}
\makeatletter\let\c@equation\c@theorem\makeatother
\newenvironment{acknowledgement}{\section*{Acknowledgement}}{}
\address{Department of Mathematics, University of Iowa, 14 MacLean Hall, Iowa City, IA 52242-1419 U.S.A.}
\email{mtakase@math.uiowa.edu}
\date{}
\author{Masamichi Takase}
\title{Homology $3$-spheres in codimension three}
\begin{document}\sloppy
\maketitle
\begin{abstract}
For smooth embeddings of an integral homology $3$-sphere in the $6$-sphere, 
we define an integer invariant in terms of their Seifert surfaces. 
Our invariant gives a bijection between the set of 
smooth isotopy classes of such embeddings and the integers. 
It also gives rise to a complete invariant 
for homology bordism classes of 
all embeddings of homology $3$-spheres in the $6$-sphere. 
As a consequence, we show that 
two embeddings of an oriented integral homology 
$3$-sphere in the $6$-sphere are isotopic if and only if 
they are homology bordant. 
We also relate our invariant to the Rohlin invariant 
and accordingly characterise those embeddings 
which are compressible into the $5$-sphere. 
\end{abstract}

\section{Introduction}
Haefliger's seminal work \cite{hae1} shows that 
in the \textit{smooth} category, 
unlike in the piecewise linear and topological categories, 
knotting phenomena occur in codimension greater than three. 
For example, \cite{hae1,hae2} show that 
the group $C^{2k+1}_{4k-1}$ of smooth isotopy classes of 
smooth embeddings of the $(4k-1)$-sphere $S^{4k-1}$ in the $6k$-sphere $S^{6k}$ 
forms the infinite cyclic group for each $k\ge1$. 
According to \cite{boechat,b-h,takase2}, 
the isotopy class of an embedding in $C^{2k+1}_{4k-1}$ 
can be read off from geometric characteristics of 
its Seifert surface. 

The first occasion $C^3_3$, where 
smooth and piecewise linear isotopies differ, 
seems especially intriguing, 
since it lies in the interface 
between high- and low-dimensional topology.  
For instance, Montgomery and Yang \cite{donaldson} showed 
an isomorphic correspondence between the group $C^3_3$ and 
the group of diffeomorphism classes of all homotopy ${\C}P^3$, 
that has provided various interesting applications (e.g. \cite{levine,lu,masuda}).   
Bo\'echat and Haefliger \cite{b-h} used the group $C^3_3$ 
to describe a necessary and sufficient condition 
for the embeddability of an orientable $4$-manifold in $7$-space 
(later their condition turned out to be always satisfied due to \cite{donaldson}, 
so that this last case of the hard Whitney embedding theorems was settled). 
Our previous paper \cite{takase3} also shows that the group $C^3_3$
is linked to various aspects of $4$-dimensional topology. 

In this paper, 
we study smooth embeddings of oriented integral homology $3$-spheres in $S^6$. 
We first observe  
in Proposition~\ref{prop:seifert} that 
for every smooth embedding $F\co\Sigma^3\hookrightarrow S^6$ of 
an oriented integral homology $3$-sphere in $S^6$, 
there exists a smooth embedding (a Seifert surface)
$\widetilde{F}\co V^4\hookrightarrow S^6$ of 
a smooth oriented $4$-manifold $V^4$ with $\partial{V^4}=\Sigma^3$ 
such that $\widetilde{F}|_{\partial{V^4}}=F$. 
For such a Seifert surface, 
we denote the signature of $V^4$ by $\sigma(V^4)$ and 
the normal Euler class of $\widetilde{F}$ by $e_{\widetilde{F}}$. 
In \S\ref{sect:invariant}, we prove that 
\[
\Omega(F)=\frac{\sigma(V^4)-e_{\widetilde{F}}\smile e_{\widetilde{F}}}{8}, 
\]
does not depend on the choice of $\widetilde{F}$ and 
is invariant up to isotopy of $F$, 
where $e_{\widetilde{F}}\smile e_{\widetilde{F}}$ is the 
square of $e_{\widetilde{F}}$, evaluated on the fundamental class. 
Then, 
Theorem~\ref{thm:isotopic} further shows  
that the invariant $\Omega$ gives a bijection: 
\[\begin{array}{ccc}
\Omega\co\Emb(\Sigma^3,S^6)\stackrel{\approx}{\to}\Z
\end{array}\]
between the set $\Emb(\Sigma^3,S^6)$ of smooth isotopy classes of 
smooth embeddings of $\Sigma^3$ in $S^6$ and the integers $\Z$. 

In the study of oriented integral homology $3$-spheres, 
their bounding \textit{smooth} compact oriented $4$-manifolds 
(note that they always bound \textit{topological} compact 
oriented acyclic $4$-manifolds \cite{freedman})
are important tools and give rise to 
many useful invariants in low-dimensional topology. 
In Corollary~\ref{cor:compress1}, 
we relate our invariant $\Omega$ to the Rohlin invariant 
and consequently show that an embedding $F$ of 
an integral homology sphere $\Sigma^3$ in $S^6$ 
is compressible into $S^5$ if and only if 
$\Omega(F)$ is modulo $2$ equal to the Rohlin invariant 
$\mu(\Sigma^3)$. 

In \S\ref{sect:cobordism} we show that 
when we assemble all embeddings of 
oriented integral homology $3$-spheres in $S^6$, 
the invariant $\Omega$ proves to be invariant up to 
homology bordism (of embeddings). 
We say that two embeddings $F_0$ and $F_1$ of two integral homology spheres 
are \textit{homology bordant} if there exists 
a proper embedding of a homology cobordism between the two homology spheres 
in $S^6\times[0,1]$, 
whose restriction to the boundary coincides with 
the disjoint union of $F_0\times\{0\}$ and $-F_1\times\{1\}$ (see \S\ref{sect:cobordism} for details). 
The collection of homology bordism classes forms 
an abelian group, denoted by $\Gamma^3_3$, via connected sum. 
Then, we show that the invariant $\Omega$ gives rise to 
the following isomorphism (Theorem~\ref{thm:cobordism}): 
\[\begin{array}{cccc}
\bar{\Omega}\co&\Gamma^3_3&\stackrel{\approx}{\to}&\Theta^3_{\Z}\,\oplus\,\Z\\
&{[F\co\Sigma^3\hookrightarrow S^6]}&\mapsto&(\,[\Sigma^3],\,\Omega(F)\,) 
\end{array}\]
where 
$[\Sigma^3]$ denotes the homology cobordism class represented 
by the homology sphere $\Sigma^3$ and 
$[F\co\Sigma^3\hookrightarrow S^6]$ denotes the homology bordism class 
represented by the embedding $F\co\Sigma^3\hookrightarrow S^6$. 
As a corollary, 
we see that two embeddings of an oriented integral homology 
$3$-sphere in $S^6$ are isotopic if and only if 
they are homology bordant. 

We work in the smooth category; 
all manifolds and mappings 
are supposed to be differentiable of class $C^\infty$, 
unless otherwise explicitly stated. 
We will suppose all spheres and homology spheres are oriented. 
If $M$ is an oriented manifold with non-empty boundary, 
then for the \textit{induced orientation of\/ $\partial{M}$} 
we adopt the \textit{outward vector first} convention: 
we say an ordered basis of 
$T_p(\partial{M})$ ($p\in\partial{M}$) is positively oriented 
if an outward vector followed by the basis
is a positively oriented basis of $T_pM$. 
For a closed $n$-dimensional manifold $M$ we denote its punctured manifold 
by $M_\circ$; i.e., $M_\circ:= M\smallsetminus\Int{D^n}$.
The homology and cohomology groups are supposed to be 
with integer coefficients unless otherwise explicitly noted. 

\section{Seifert surfaces}
The purpose of this section is to show that 
every embedding $F\co\Sigma^3\hookrightarrow S^6$ of 
an oriented integral homology $3$-sphere $\Sigma^3$ has a Seifert surface. 

\begin{definition}
Let $F\co M^3\hookrightarrow S^6$ be an embedding of a closed oriented $3$-manifold. 
Then, \textit{a Seifert surface for $F$} is 
an embedding $\widetilde{F}\co W^4\hookrightarrow S^6$ 
of a compact connected oriented $4$-manifold $W^4$ with $\partial{W^4}=M^3$ 
such that $\widetilde{F}|_{M^3}=F$. 
\end{definition}

Every closed oriented $3$-manifold bounds a spin $4$-manifold 
with only one $0$-handle and $2$-handles and 
such a $4$-manifold can be embedded in $S^6$, in fact, in $S^5$ 
(e.g. see \cite[Chapter VII]{kirby}). 
Now let $\Sigma^3$ be an oriented integral homology $3$-sphere and 
take an embedding $\widetilde{F}_0\co W^4_0\hookrightarrow S^6$ 
of a compact oriented $4$-manifold $W^4_0$ 
with $\partial{W^4_0}=\Sigma^3$. 
In this section, we fix this embedding $\widetilde{F}_0$ and 
put $F_0:=\partial{\widetilde{F}_0}\co\Sigma^3\hookrightarrow S^6$. 
We can assume that $F_0$ is standard on a $3$-disc $D^3\subset M^3$ 
(i.e., coincides with the ``northern hemisphere'' $D^3_+\subset S^6_+$ 
of the standard inclusion $S^3\subset S^6$) and 
$F_0(\Sigma^3\smallsetminus\Int{D^3})$ lies in 
the southern hemisphere $D^6_-$ of $S^6$. 
Note that the complement $C$ of an open neighbourhood of $F_0(\Sigma^3\smallsetminus\Int{D^3})$ 
in $S^6$ is contractible. 

First, we need the following lemma. 

\begin{lemma}\label{lem:punc}
Let $M^3$ be a closed oriented $3$-manifold. 
Then, two embeddings $M^3_\circ\hookrightarrow S^6$ 
of the punctured manifold $M^3_\circ$ in $S^6$ are isotopic. 
\end{lemma}

\begin{proof}
The punctured manifold $M^3_\circ$ has a handle decomposition 
\[
M^3_\circ=D^3\cup\{\text{$1$-handles}\}\cup\{\text{$2$-handles}\}.
\]
For dimensional reasons, two embeddings $M^3_\circ\hookrightarrow S^6$ 
are isotopic on the cores of the handles and 
we can extend the isotopy on the entire handles 
since $\pi_1(V_{5,2})=0$ and $\pi_2(V_{4,1})=0$, 
where $V_{n,k}$ denotes the Stiefel manifold of all $k$-frames in $n$-space. 
\end{proof}

\begin{remark}
If $M^3=\Sigma^3$ is an integral homology sphere, then an embedding 
$\Sigma^3_\circ\hookrightarrow S^6$ of the punctured manifold 
can be compressed in $S^5$ \cite[Theorem~D]{hirsch2} 
and is unique up to isotopy in $S^5$ \cite[Corollary~4.10(1)]{saeki}. 
\end{remark}

Next, we prove the following. 

\begin{figure}[tbp]
\begin{center}
\includegraphics[width=.75\textwidth,keepaspectratio,bb=86 625 550 808]{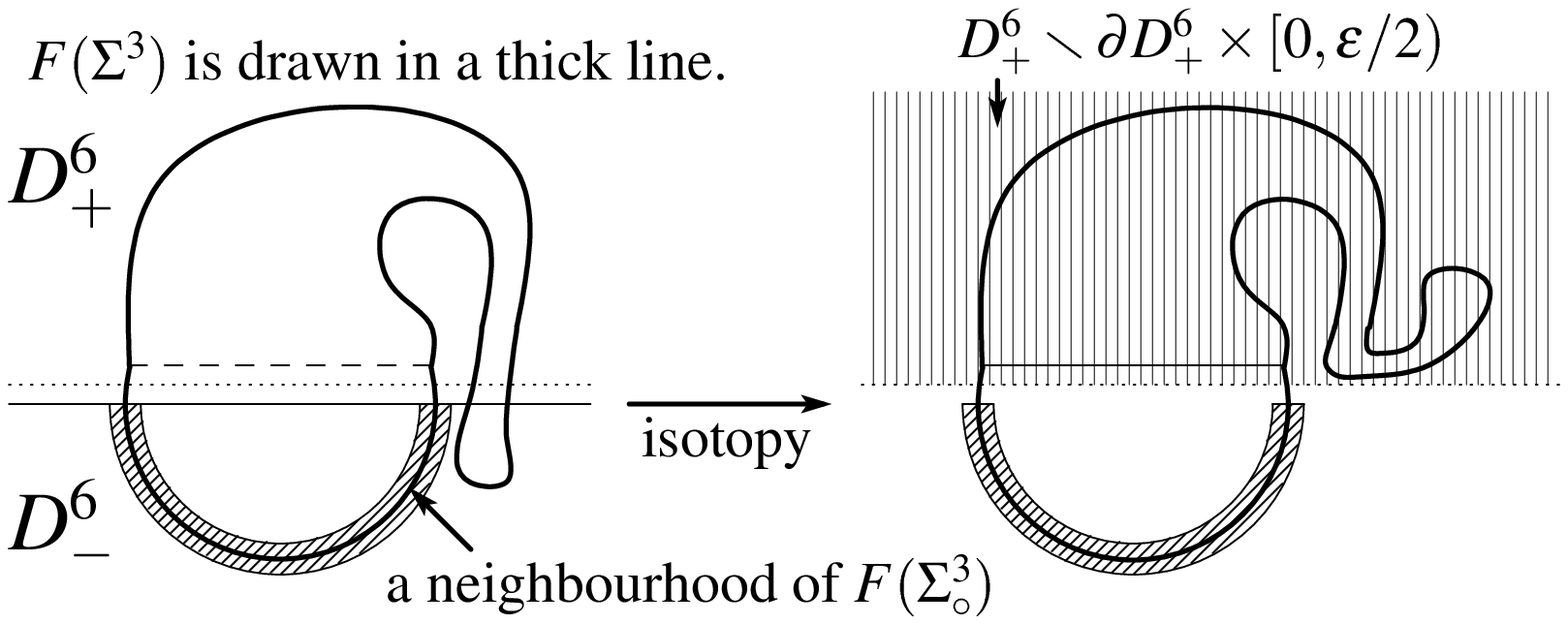}
\end{center}
\caption{}\label{fig:engulf}
\end{figure}

\begin{proposition}\label{prop:conn}
Any embedding $F\co\Sigma^3\hookrightarrow S^6$ of 
an oriented integral homology $3$-sphere is isotopic to 
the connected sum $F_0\sharp g$ of $F_0$ with some embedding 
$g\co S^3\hookrightarrow S^6$ of the $3$-sphere $S^3$. 
\end{proposition}

\begin{proof}
By Lemma~\ref{lem:punc}, $F\co\Sigma^3\hookrightarrow S^6$ 
can be isotoped to $F_0$ on 
${\Sigma^3_\circ}:=\Sigma^3\smallsetminus\Int{D^3}$. 
Then, we can assume that 
a small collar neighbourhood $\partial{D^3}\times[0,\epsilon]$ 
of $\partial{D^3}\subset D^3$ 
is ``standardly embedded'' 
in a small collar neighbourhood $\partial{D^6_+}\times[0,\epsilon]$ 
of the northern hemisphere $D^6_+$ of $S^6$ 
(i.e.,  $F$ on $\partial{D^3}\times[0,\epsilon]$ is the product of 
the standard inclusion $S^2\subset S^5=\partial{D^6_+}$ by the interval). 

The complement $C$ of an open neighbourhood of $F(\Sigma^3_\circ)$ 
in $S^6$ is contractible. 
Therefore, we can apply a lemma by Levine \cite[Lemma~2]{levine} 
to $C$ and its two submanifolds, 
$D^6_+\smallsetminus\partial{D^6_+}\times[0,\epsilon/2)$ 
in codimension $0$ and 
$F(D^3\smallsetminus\partial{D^3}\times[0,\epsilon))$ 
capped off by a $3$-disc in $\partial{D^6_+}\times\{\epsilon\}$ 
(with corners smoothed) in codimension $3$ (see Figure~\ref{fig:engulf}). 
This ensures that we can push $F(D^3)$ into $D^6_+\subset C$ by 
an isotopy of $C$ fixing a collar of $F(D^3)$. 
Hence, $F$ is isotopic to 
the connected sum $F_0\conn g$ of 
$F_0$ with some embedding $g\co S^3\hookrightarrow S^6$. 
\end{proof}

The following is now an easy corollary of Proposition~\ref{prop:conn}. 

\begin{proposition}\label{prop:seifert}
Any embedding $F$ of an oriented integral homology $3$-sphere $\Sigma^3$ in $S^6$ 
has a Seifert surface. 
\end{proposition}

\begin{proof}
By Proposition~\ref{prop:conn}, 
$F$ is isotopic to 
the connected sum $F_0\conn g$ of 
$F_0$ with some embedding $g\co S^3\hookrightarrow S^6$. 
By \cite{g-m,takase2}, $g$ has a Seifert surface 
$\widetilde{g}\co V^4\hookrightarrow S^6$. 
Thus, the boundary connected sum 
$\widetilde{F}_0\bconn \widetilde{g}\co W^4_0\bconn V^4\hookrightarrow S^6$ 
(composed with a suitable diffeomorphism on $S^6$) 
is a Seifert surface for $F=F_0\conn g$. 
\end{proof}

\begin{remark}\label{rem:seifert} 
In Proposition~\ref{prop:seifert}, 
we can choose both $W_0$ and $V^4$ to be simply connected 
(in fact, $V^4$ can be chosen as $S^2\times S^2$ \cite{takase2} or 
as the connected sum of some copies of the complex projective planes \cite{g-m,takase3}). 
Therefore, we see that 
any embedding of an oriented integral homology $3$-sphere 
has a simply connected Seifert surface. 
\end{remark}

\section{The Hopf invariant}
In this section 
we associate to each embedding of a homology $3$-sphere into $S^6$, 
which is equipped with a normal vector field, its \textit{Hopf invariant}.

\subsection{Embeddings with normal $1$-fields}\label{subsect:1-field} 
Let $F\co M^3\hookrightarrow S^6$ be an embedding of 
a closed oriented $3$-manifold $M^3$. 
Then, the normal bundle of $F$ is trivial and 
homotopy classes of normal framings 
(i.e.\ trivialisations of the normal bundle) are classified by 
the homotopy set $[M^3,\SO(3)]$, with respect to some fixed framing. 
Similarly, homotopy classes of normal $1$-fields 
are classified by the set $[M^3,V_{3,1}]=[M^3,S^2]$, 
which has been computed by Pontrjagin \cite{pontrjagin}. 

If $\Sigma^3$ is an integral homology sphere, 
then the absence of $2$-dimensional obstructions allows us to 
identify $[\Sigma^3,S^2]$ with $H^3(\Sigma^3)=\Z$ in a natural way. 
In this case furthermore, 
the natural map $[\Sigma^3,\SO(3)]\to[\Sigma^3,S^2]$ 
is a bijection (see \cite[Proposition~2.2]{kuperberg}). 
Therefore, a normal $1$-field determines 
a normal framing uniquely up to homotopy and vice versa. 
Hence, in a similar way to \cite[\S5.1]{hae2}, 
we can consider the connected sum of two 
embeddings $\Sigma^3\hookrightarrow S^6$ with normal $1$-fields. 

\subsection{The Hopf invariant of an embedding with a normal $1$-field}\label{subsect:hopf}
Let $\Sigma^3$ be an oriented integral homology $3$-sphere 
and $F\co \Sigma^3\hookrightarrow S^6$ be an embedding. 
Then, by Alexander duality, 
the exterior space $X:=S^6\smallsetminus F(\Sigma^3)$ has 
the same homology groups as the $2$-sphere $S^2$.  
We identify $H_2(X)$ with $\Z$ by assigning to 
$[z^2]\in H_2(X)$ the linking number 
$\lk(z^2,f(\Sigma^3))\in\Z$ in $S^6$, 
which also gives an identification $H^2(X)=\Z$. 

If an embedding $F\co \Sigma^3\hookrightarrow S^6$ is 
endowed with a normal $1$-field $\nu$ of $F(\Sigma^3)\subset S^6$, 
we can define the Hopf invariant of $(F,\nu)$ in the following way. 

\begin{definition}\label{defn:fcp}
Let $F\co \Sigma^3\hookrightarrow S^6$ be an embedding 
of an oriented homology $3$-sphere and 
$\nu$ be a normal $1$-field of $F(\Sigma^3)\subset S^6$. 
By using a small shift of $F(\Sigma^3)$ along $\nu$, 
we can define a map $\bar{\nu}\co\Sigma^3\to X:=S^6\smallsetminus F(\Sigma^3)$. 
Then, we define \textit{the Hopf invariant $H_{(F,\nu)}$ of $(F,\nu)$} to be minus 
\textit{the functional cup product} 
$L_{\bar{\nu}}(\omega,\omega)=R_{\bar{\nu}}(\omega,\omega)\in H^3(\Sigma^3)=\Z$ 
for the generator $\omega\in H^2(X)$. 
We use the definition of the functional cup product 
given in \cite[p.368]{u-m} and \cite[\S2.2]{takase2}. 
In our special case here, we can review it as follows. 

Take a $2$-cochain $\omega'$ of $X$ which represents $\omega\in H^2(X)$. 
Then, there exist 
a $3$-cochain $a$ of $X$ such that $\delta{a}=\omega'\smile\omega'$ and 
a $1$-cochain $b$ of $\Sigma^3$ such that $\delta{b}=\bar{\nu}^{\sharp}\omega'$, 
where $\bar{\nu}^{\sharp}$ is the induced homomorphism on cochains. 
Finally, the $3$-cochain 
\begin{equation}\label{eqn:fcp}
z'=\bar{\nu}^{\sharp}a-b\smile\bar{\nu}^{\sharp}\omega'
\end{equation}
actually becomes a cocycle of $\Sigma^3$. 
Thus, we define the functional cup product 
$L_{\bar{\nu}}(\omega,\omega)$ to be $[z']\in H^3(\Sigma^3)$ 
and the Hopf invariant $H_{(F,\nu)}$ of $(F,\nu)$ to be 
$-[z']\in H^3(\Sigma^3)=\Z$, 
where we identify $H^3(\Sigma^3)$ with the integers $\Z$ 
with respect to the orientation of $\Sigma^3$. 
Note that $H_{(F,\nu)}$ does not depend on the choice of 
the orientation of $X$. 
\end{definition}

\begin{definition}
Let $F\co \Sigma^3\hookrightarrow S^6$ be an embedding 
of an oriented integral homology $3$-sphere and 
$\widetilde{F}\co V^4\hookrightarrow S^6$ be its Seifert surface. 
Then, we define \textit{the Hopf invariant $H_{\widetilde{F}}$ for $\widetilde{F}$} 
to be the Hopf invariant of 
\[(F,\text{the outward normal field of $F(\Sigma^3)\subset\widetilde{F}(V^4)$}).\]
\end{definition}

\begin{theorem}\label{thm:square}
Let $F\co \Sigma^3\hookrightarrow S^6$ be an embedding 
of an oriented integral homology $3$-sphere and 
$\widetilde{F}\co V^4\hookrightarrow S^6$ be its Seifert surface. 
Then, 
\[H_{\widetilde{F}}=-e_{\widetilde{F}}\smile e_{\widetilde{F}}\in\Z,\]
where $e_{\widetilde{F}}\in H^2(V^4)=H^2(V^4,\partial{V^4})$ is 
the normal Euler class of $\widetilde{F}$ and 
$e_{\widetilde{F}}\smile e_{\widetilde{F}}$ is its square 
evaluated on the fundamental homology class of $V^4$. 
\end{theorem}

\begin{proof}
We calculate the Hopf invariant by translating 
the cup product operation on cohomology 
into intersection theory on homology 
using the duality theorems. 

Put $\HS:=F(\Sigma^3)$ and denote by $\HS'\in S^6\smallsetminus\HS$ 
a small shift of $\HS$ along the outward normal field of 
$F(\Sigma^3)\subset\widetilde{F}(V^4)$. 
Let ${\bar{\nu}}\co\HS'\hookrightarrow S^{6}\smallsetminus\HS$ 
be the inclusion and 
\[\omega=[\omega']\in H^{2}(S^{6}\smallsetminus\HS)=\Z\]
denote the generator (see \S\ref{subsect:hopf}). 

Let $N$, $N'$ be sufficiently small (i.e.\,$N\cap N'=\emptyset$) 
tubular neighbourhoods of 
$\HS$, $\HS'$ respectively, and
put $X:=S^{6}\smallsetminus\Int{N}$. 
Then, consider the following diagram
\[\begin{CD}
H^*(S^{6}\smallsetminus\HS)@>{\bar{\nu}}^*>>H^*(\HS')\\
@V{\approx}VV @AA{\approx}A\\
H^*(X)@>{\hat{\nu}}^*>>H^*(N'),
\end{CD}\]
where all homomorphisms are induced by the inclusion maps and 
the vertical arrows are isomorphisms. 
By this diagram, we can calculate the functional cup product 
\[L_{\hat{\nu}}(\omega,\omega)=R_{\hat{\nu}}(\omega,\omega)\in H^{3}(N')=\Z\]
with respect to the inclusion ${\hat{\nu}}\co N'\hookrightarrow X$, 
instead of $L_{\bar{\nu}}(\omega,\omega)=R_{\bar{\nu}}(\omega,\omega)\in H^{3}(\HS')=\Z$, 
where we identify $\omega$ with its inverse image 
(also denoted by $\omega$) under the isomorphism 
$H^{2}(S^{6}\smallsetminus\HS)\to H^{2}(X)$.

Furthermore, by the duality and the excision theorems, 
we have the following:
\[\begin{CD}
H^i(X)@>{\hat{\nu}}^*>>H^i(N')\\
@V{\approx}VV @VV{\approx}V\\
H_{6-i}(X,\partial{X})@>{\hat{\nu}}_!>>H_{6-i}(N',\partial N')\\
@A{\approx}AA \\
H_{6-i}(S^{6},\HS),
\end{CD}\]
where all vertical arrows are isomorphisms.
Now using the above diagram, 
we calculate the desired Hopf invariant 
in terms of intersections on homology groups. 

The class in 
$H_4(X,\partial X)\approx H_4(S^{6},\HS)$ 
dual to the generator 
$\omega\in H^{2}(X)$ 
is represented by $\mathtt{V}:=\widetilde{F}(V^{4})$. 
Extend the outward normal field of $F(\Sigma^3)\subset\widetilde{F}(V^4)$ 
to a vector field in $S^6$ and perturb $\mathtt{V}$ by this field 
into $\mathtt{V}'$. 
We can assume that the intersection 
$\Delta:=\mathtt{V}\cap\mathtt{V}'$ of $\mathtt{V}$ and $\mathtt{V}'$ 
lies in their interior; 
$\Delta=\Int{\mathtt{V}}\cap\Int{\mathtt{V}'}$ 
(since the normal bundle of $F$ is trivial). 
Then, $\widetilde{F}^{-1}(\Delta)\subset V^4$ 
represents the integral 
dual to the normal Euler class $e_{\widetilde{F}}$ of $\widetilde{F}$. 

Since a sufficiently small tubular neighbourhood $N'$ of 
$\HS'\subset S^{6}$ does not intersect $\mathtt{V}$, 
the relative $4$-chain of $(N',\partial{N'})$
corresponding to ${\bar{\nu}}^\sharp\omega'$ in 
the second term of the right-hand side of (\ref{eqn:fcp}) 
in Definition~\ref{defn:fcp} vanishes. 

Next we check the first term of the right-hand side in (\ref{eqn:fcp}). 
Since the relative $2$-chain dual to $\omega'\smile\omega'$ 
is represented by the intersection $\Delta$ ($=\mathtt{V}\cap\mathtt{V}'$), 
the term corresponding to $\bar{\nu}^{\sharp}a$ (in (\ref{eqn:fcp})) 
is represented by the intersection of 
$\HS$ and a relative $3$-chain bounded by $\Delta$. 
Therefore, the homology class dual to 
the functional cup product $L_{\hat{\nu}}(\omega,\omega)$ is equal to 
\[
-\text{(a relative $3$-chain bounded by $\Delta$)}\cap\HS=-\lk(\Delta,\HS)\\
=-[\Delta]\in H_{2}(X)=\Z, 
\]
where the minus sign is due to the change of the order of the intersection. 
Thus, we obtain $H_{\widetilde{F}}=[\Delta]$. 

If we take another copy $\mathtt{V}''$ 
of $\mathtt{V}=\widetilde{F}(V^{4k})$ 
perturbed in an appropriate manner 
and put $\Delta':=\mathtt{V}\cap\mathtt{V}''$, 
then we have 
\[\begin{array}{rll}
[\Delta]
&=[\Delta'] &(\in H_{2}(X)=\Z)\\
&=\lk(\Delta',\HS) &(\in H_0(S^{6k})=\Z)\\
&=-\lk(\HS,\Delta')\\
&=-[\mathtt{V}\cap\Delta']\\
&=-[\mathtt{V}\cap\mathtt{V}'\cap\mathtt{V}'']\\
&=-[\Delta\cap\Delta'] &(\in H_0(\mathtt{V})=H_0(V^4)=\Z)\\
&=-e_{\widetilde{F}}\smile e_{\widetilde{F}} &(\in H^{4}(\widehat{V}^{4})=\Z).
\end{array}\]
This completes the proof.
\end{proof}

\begin{remark}
Theorem~\ref{thm:square} is very similar to \cite[Theorem~5.1]{takase2}, 
which however contained a sign error. 
To be precise, the error occurs in the last sentence of the proof of \cite[Lemma~5.2]{takase2}. 
It reads ``\textit{Thus, we have 
$H_{\widetilde{F}}=-[\widetilde{F}(\Sigma^{2k})].$}'' but 
should be 
``\textit{Thus, the desired functional product is computed to be equal 
to $-[\widetilde{F}(\Sigma^{2k})]$ and hence the Hopf invariant 
$H_{\widetilde{F}}$ is equal to $[\widetilde{F}(\Sigma^{2k})]\in H_{2k}(X)=\Z$}.'' 
With this correction, each term involving the square of the normal Euler class, 
appearing in \cite[Theorem~5.1, Corrollaries~6.2, 6.3(a) and 6.5]{takase2}, 
should change its sign. 
\end{remark}

\subsection{The case of embeddings of the $3$-sphere}
Haefliger \cite{hae1} and \cite[\S5.16]{hae2} 
proved that the group $C^3_3$ of isotopy classes of embeddings 
$S^3\hookrightarrow S^6$ is isomorphic to the integers $\Z$. 
In fact, he gave a complete invariant $\Omega\co C^3_3\to\Z$ 
and also an explicit construction of an embedding 
representing a generator of $C^3_3$. 

The following is given in 
\cite{g-m,takase2} (see also \cite{b-h,takase3})

\begin{theorem}\label{thm:3dim}
Every embedding $F\co S^3\hookrightarrow S^6$ has a Seifert surface 
$\widetilde{F}\co V^4\hookrightarrow S^6$ and 
\begin{eqnarray*}
\Omega(F)
&=&-\frac{1}{8}(\sigma(V^4)+H_{\widetilde{F}})\\
&=&-\frac{1}{8}(\sigma(V^4)-e_{\widetilde{F}}\smile e_{\widetilde{F}})
\end{eqnarray*} 
gives the isomorphism $\Omega\co C^3_3\stackrel{\approx}{\to}\Z.$ 
\end{theorem}

The following two examples of Seifert surfaces for the generator in $C^3_3$, 
given in \cite[\S2.3]{takase3}, 
will be important in our argument. 
Let ${\C{P}^2}_{\circ}$ be the punctured complex projective plane 
and $\overline{\C{P}^2}_{\circ}$ 
the punctured complex projective plane with the reversed orientation. 

\begin{proposition}\label{prop:p}
The standard inclusion $S^3\subset S^6$ has a Seifert surface 
$\widetilde{P}\co{\C{P}^2}_{\circ}\hookrightarrow S^6$ 
with Hopf invariant $H_{\widetilde{P}}=-1$. 
\end{proposition}

\begin{proposition}\label{prop:q}
The standard inclusion $S^3\subset S^6$ has a Seifert surface 
$\widetilde{Q}\co\overline{\C{P}^2}_{\circ}\hookrightarrow S^6$ 
with Hopf invariant $H_{\widetilde{Q}}=1$. 
\end{proposition}

\begin{remark}\label{rem:pq}
Propositions~\ref{prop:p} and \ref{prop:q} imply that 
by taking the boundary connected sums 
with a suitable number of copies of $\widetilde{P}$ or of $\widetilde{Q}$, 
we can alter a given Seifert surface for an embedding 
$F\co S^3\hookrightarrow S^6$ 
so that it has an arbitrarily preassigned Hopf invariant 
(without changing the isotopy class of $F$). 
This means further that for a given 
$F\co S^3\hookrightarrow S^6$ 
we can choose a Seifert surface so that near the boundary 
it extends $F$ in an arbitrarily prescribed normal direction, 
since the natural homomorphism $\pi_3(SO(3))\to\pi_3(S^2)$ 
and the Hopf invariant $\pi_3(S^2)\to\Z$ 
are both isomorphisms (see \cite{takase3}). 
\end{remark}

\subsection{The additivity of the Hopf invariant}
\begin{figure}[tbp]
\begin{center}
\includegraphics[width=.75\textwidth,keepaspectratio,bb=119 668 425 715]{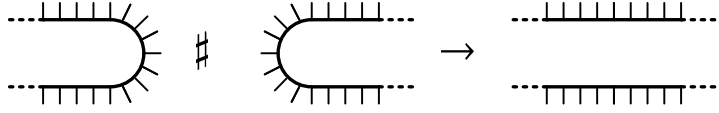}
\end{center}
\caption{The connected sum of embeddings with normal $1$-fields}\label{fig:conn}
\end{figure}
Let $(F_i,\nu_i)\co \Sigma^3_i\hookrightarrow S^6$ ($i=1,2$) be embeddings 
with normal $1$-fields of homology $3$-spheres. 
Then, we can consider the connected sum $(F_1,\nu_1)\conn(F_2,\nu_2)$
in a natural way 
(see Figure~\ref{fig:conn}) and we have the following. 

\begin{proposition}\label{prop:additivity}
Let $(F_i,\nu_i)\co \Sigma^3_i\hookrightarrow S^6$ $(i=1,2)$ be 
two embeddings with normal $1$-fields of homology $3$-spheres. Then, 
\[
H_{(F_1,\nu_1)\conn(F_2,\nu_2)}=H_{(F_1,\nu_1)}+H_{(F_2,\nu_2)}\quad\in\Z. 
\]
\end{proposition}

To prove Proposition~\ref{prop:additivity}, we use the following lemma, 
which will be also needed later.  

\begin{lemma}\label{lem:normal}
Let $F\co \Sigma^3\hookrightarrow S^6$ be an embedding 
of a homology $3$-sphere $\Sigma^3$ and 
$\nu$ be a normal $1$-field of $F(\Sigma^3)\subset S^6$. 
Then, there exists a Seifert surface 
$\widetilde{F}\co W^4\hookrightarrow S^6$ 
for $F$ such that 
the outward normal field of $F(\Sigma^3)\subset\widetilde{F}(W^4)$ 
coincides with $\nu$. 
\end{lemma}

\begin{proof}
Let $\widetilde{F}'\co W'\hookrightarrow S^6$
be a Seifert surface for $F$ by Proposition~\ref{prop:seifert} and 
denote by $\nu'$ 
the outward normal vector field of $F(\Sigma^3)\subset\widetilde{F}'(W')$. 
Then, $\nu'$ and the given $\nu$ are homotopic 
on ${\Sigma^3_\circ}=\Sigma^3\smallsetminus\Int{D^3}$, 
since $H^2(\Sigma^3)=0$ (cf.\ \cite[Corollary~4.9]{gompf}). 
We can resolve the difference between $\nu'$ and $\nu$  on the final $3$-cell 
by taking the boundary connected sum of $\widetilde{F}'$ 
with a suitable number of copies of $\widetilde{P}$ or of $\widetilde{Q}$
(in Propositions~\ref{prop:p} and \ref{prop:q}). 
Thus, we have a new Seifert surface for $F$ 
for which the outward vector field along the boundary 
is homotopic to $\nu$. 
After composing a diffeomorphism of $S^6$ if necessary, 
we obtain a desired Seifert surface. 
\end{proof}

\begin{remark}\label{rem:normal}
In view of Remark~\ref{rem:seifert}, 
the Seifert surface in Lemma~\ref{lem:normal} 
can always chosen to be simply connected. 
\end{remark}

Now Proposition~\ref{prop:additivity} is an easy corollary of 
Theorem~\ref{thm:square} and Lemma~\ref{lem:normal}. 

\begin{proof}[Proof of Proposition~\ref{prop:additivity}]
By Lemma~\ref{lem:normal}, 
for each $(F_i,\nu_i)\co \Sigma^3_i\hookrightarrow S^6$ 
we can consider a Seifert surface 
$\widetilde{F}_i\co W^4_i\hookrightarrow S^6$ 
so that the outward normal vector field of 
$F_i(\Sigma_i)\subset\widetilde{F}_i(W^4_i)$
coincides with $\nu_i$. 

Then, by Theorem~\ref{thm:square}, 
\[
H_{(F_1,\nu_1)\conn(F_2,\nu_2)}=-e_{\widetilde{F}_1\bconn\widetilde{F}_2}\smile e_{\widetilde{F}_1\bconn\widetilde{F}_2} 
\]
where $e_{\widetilde{F}_1\bconn\widetilde{F}_2}$ is the normal Euler class of $\widetilde{F}_1\bconn\widetilde{F}_2$. 
The homology class dual to $e_{\widetilde{F}_1\bconn\widetilde{F}_2}$ is represented by 
the intersection $\Delta$ of $(\widetilde{F}_1\bconn\widetilde{F}_2)(W^4_1\bconn W^4_2)$ 
with its small perturbation in $S^6$. 
Since $\widetilde{F}_1\bconn\widetilde{F}_2$ is an embedding and 
$F_1\conn F_2$ has trivial normal bundle, 
$\Delta$ can be assumed to be the disjoint union 
$\Delta_1\sqcup\Delta_2$ of $\Delta_i\subset\widetilde{F}_i(\Int{W^4_i})$ 
dual to $e_{\widetilde{F}_i}$ ($i=1,2$). 
Thus, 
\begin{eqnarray*}
H_{(F_1,\nu_1)\conn(F_2,\nu_2)}
&=&-[\Delta]\bullet[\Delta]\\
&=&-[\Delta_1]\bullet[\Delta_1]-[\Delta_2]\bullet[\Delta_2]\\
&=&-e_{\widetilde{F}_1}\smile e_{\widetilde{F}_1}-e_{\widetilde{F}_2}\smile e_{\widetilde{F}_2}\\
&=&H_{\widetilde{F}_1}+H_{\widetilde{F}_2}\\
&=&H_{(F_1,\nu_1)}+H_{(F_2,\nu_2)}, 
\end{eqnarray*}
where $\bullet$ means the intersection pairing of homology. 
\end{proof}

Furthermore, we have the following. 

\begin{corollary}\label{cor:inj}
Let $F\co \Sigma^3\hookrightarrow S^6$ be an embedding 
of a homology $3$-sphere $\Sigma^3$. Then, 
two normal $1$-fields $\nu_1$ and $\nu_2$ of $F(\Sigma^3)\subset S^6$ 
are homotopic if and only if 
$H_{(F,\nu_1)}=H_{(F,\nu_2)}$. 
\end{corollary}

\begin{proof}
Since $H^2(\Sigma^3;\pi_2(S^2))=0$, 
the $1$-fields 
$\nu_1$ and $\nu_2$ are homotopic on $\Sigma^3\smallsetminus\Int{D^3}$. 
Therefore, 
their first difference may be in 
$H^3(\Sigma^3;\pi_3(S^2))=\Z$ and 
we can assume that $(F,\nu_2)=(F,\nu_1)\conn(j,\nu)$ 
for the standard inclusion $j\co S^3\hookrightarrow S^6$ 
and a normal $1$-field $\nu$ of $j(S^3)\subset S^6$. 
Then, $\nu_1$ and $\nu_2$ are homotopic if and only if 
$\nu$ is \textit{standard}, i.e., is homotopic to 
the first vector field to the standard normal framing 
of the standard inclusion $j$. 
Since the Hopf invariant $\pi_3(S^2)\to\Z$ is an isomorphism, 
we have further 
\begin{eqnarray*}
\text{$\nu_1$ and $\nu_2$ are homotopic}
&\Leftrightarrow&\text{$\nu$ is standard}\\
&\Leftrightarrow&H_{(j,\nu)}=0\\
&\Leftrightarrow&H_{(F,\nu_2)}=H_{(F,\nu_1)}+H_{(j,\nu)}=H_{(F,\nu_1)}. 
\end{eqnarray*}
This completes the proof. 
\end{proof}

\section{The invariant}\label{sect:invariant}
In the section, 
we prove the analogue of Theorem~\ref{thm:3dim} 
for embeddings of an integral homology sphere in $S^6$. 

First, we prove the following. 

\begin{proposition}\label{prop:inv}
Let $F\co\Sigma^3\hookrightarrow S^6$ be an embedding of 
an oriented integral homology $3$-sphere and 
$\widetilde{F}\co W^4\hookrightarrow S^6$ 
be its Seifert surface. 
Then, the integer 
\begin{eqnarray*}
\Omega(F)
&:=&-\frac{1}{8}(\sigma(W^4)+H_{\widetilde{F}})\\
&=&-\frac{1}{8}(\sigma(W^4)-e_{\widetilde{F}}\smile e_{\widetilde{F}}). 
\end{eqnarray*} 
does not depend on the choice of $\widetilde{F}$ and 
depends only on the isotopy class of $F$. 
\end{proposition}

\begin{proof}
First, notice that 
since the modulo two reduction of $e_{\widetilde{F}}$ equals 
the second Stiefel-Whitney class $w_2(W^4)$ of $W^4$, 
the integer $\sigma(W^4)-e_{\widetilde{F}}\smile e_{\widetilde{F}}$, 
and hence $\sigma(W^4)+H_{\widetilde{F}}$ (by Theorem~\ref{thm:square}) 
are divisible by $8$. 

Consider two mutually isotopic embeddings 
$F_i\co \Sigma^3_i\hookrightarrow S^6$ 
and their Seifert surfaces 
$\widetilde{F}_i\co W^4_i\hookrightarrow S^6$ ($i=1,2$).  
Then, there exists a diffeomorphism $h$ on $S^6$ 
such that $F_1=h\circ F_2$. 
We abuse the symbols $F_2$ and $\widetilde{F}_2$ 
to denote respectively $h\circ F_2(=F_1)$ and $h\circ\widetilde{F}_2$. 

Put $k:=H_{\widetilde{F}_1}-H_{\widetilde{F}_2}$ and 
consider the boundary connected sum $\widetilde{F}_1\bconn k\widetilde{P}$, 
where $k\widetilde{P}$ means 
\[\begin{cases}
\text{the standard inclusion}&\text{if $k=0$,}\\
\text{the boundary connected sum of $k$ copies of $\widetilde{P}$}&\text{if $k>0$,}\\
\text{the boundary connected sum of $|k|$ copies of $\widetilde{Q}$}&\text{if $k<0$} 
\end{cases}\]
(see Propositions~\ref{prop:p} and \ref{prop:q}).
Then, $\widetilde{F}_1\bconn k\widetilde{P}$ is still a Seifert surface for 
$F_1$ and has the Hopf invariant 
$H_{\widetilde{F}_1\bconn k\widetilde{P}}=H_{\widetilde{F}_1}-k=H_{\widetilde{F}_2}$ 
by Proposition~\ref{prop:additivity}. 
This means that 
the outward normal fields along the boundaries 
of $\widetilde{F}_1\bconn k\widetilde{P}$ and of $\widetilde{F}_2$ 
are homotopic by Corollary~\ref{cor:inj}. 
Since a map from $\Sigma^3$ to $S^2$ and its composition with 
the antipodal map on $S^2$ are homotopic (see \S\ref{subsect:hopf}), 
Corollary~\ref{cor:inj} implies that 
we can isotope $\widetilde{F}_1\bconn k\widetilde{P}$ 
so that near the boundary it is prolonged in the outward normal direction 
of $F_2(\Sigma^3)\subset\widetilde{F}_2(W^4_2)$. 

Thus, by using $\widetilde{F}_1\bconn k\widetilde{P}$ and $\widetilde{F}_2$, 
we obtain an immersion 
\[
(\widetilde{F}_1\bconn k\widetilde{P})\cup\widetilde{F}_2\co 
(W^4_1\bconn k{\C{P}^2}_\circ)\cup_\partial(-W^4_2)\looparrowright S^6,
\]
where $(W^4_1\bconn k{\C{P}^2}_\circ)\cup_\partial(-W^4_2)$ 
denotes the closed $4$-manifold obtained by gluing 
$W^4_1\bconn k{\C{P}^2}_\circ$ and $W^4_2$ (with the reversed orientation) 
via the orientation-reversing diffeomorphism on the common boundary $\Sigma^3$. 

Clearly, this immersion $(\widetilde{F}_1\bconn k\widetilde{P})\cup\widetilde{F}_2$ 
has no triple points, since its multiple points consist only of 
the intersection of the two embeddings 
$\widetilde{F}_1\bconn k\widetilde{P}$ and $\widetilde{F}_2$.  
Hence, by \cite{white} (see also \cite[p.44]{kirby}), we have 
\begin{eqnarray*}
0
&=&\text{the algebraic number of triple points}\\
&=&\sigma((W^4_1\bconn k{\C{P}^2}_\circ)\cup_\partial(-W^4_2))\\
&=&\sigma(W^4_1)+k\sigma({\C{P}^2}_\circ)-\sigma(W^4_2)\\
&=&\sigma(W^4_1)-\sigma(W^4_2)+k\\
&=&(\sigma(W^4_1)+H_{\widetilde{F}_1})-(\sigma(W^4_2)+H_{\widetilde{F}_2}). 
\end{eqnarray*}
This implies that for an embedding $F\co\Sigma^3\hookrightarrow S^6$, 
$\Omega(F)$ does not depend on the choice of 
its Seifert surface and clearly is invariant up to isotopy. 
\end{proof}

The following Theorem~\ref{thm:isotopic} shows that 
the invariant $\Omega$ is a complete invariant up to isotopy. 
Let $\Emb(\Sigma^3,S^6)$ be the set of isotopy classes 
of embeddings of an oriented integral homology $3$-sphere $\Sigma^3$ in $S^6$. 

\begin{theorem}\label{thm:isotopic}
Let $F\co\Sigma^3\hookrightarrow S^6$ be an embedding of 
an oriented integral homology $3$-sphere and $\widetilde{F}\co V^4\hookrightarrow S^6$ 
be its Seifert surface. 
Then, 
$\Omega$ gives a bijection 
\[\Omega\co\Emb(\Sigma^3,S^6)\stackrel{\approx}{\to}\Z.\]
\end{theorem}

\begin{proof}
Since the signature and the Hopf invariant are both additive under 
the boundary connected sum of Seifert surfaces, 
the surjectivity of $\Omega$ follows from the fact that 
in the case of embeddings of the $3$-sphere $S^3$, 
the same formula gives the bijection 
$\Omega\co C^3_3\to\Z$ (see Theorem~\ref{thm:3dim}). 

We prove injectivity. 
Let $F_i\co\Sigma^3\hookrightarrow S^6$ $(i=0,1)$ be two embeddings. 
By an argument similar to that of Proposition~\ref{prop:conn}, 
$F_1$ is isotopic to $F_0\bconn g$ for some embedding $g\co S^3\hookrightarrow S^6$. 
Clearly, if $F_0$ and $F_1$ are not isotopic to each other, 
then $g$ is not isotopic to the standard embeddings
(see \cite[Lemma~1.3]{hae2}); $\Omega(g)\ne0$. 
Hence, $\Omega(F_1)=\Omega(F_0\bconn g)=\Omega(F_0)+\Omega(g)\ne\Omega(F_0)$. 
\end{proof}

\begin{remark}
Hausmann \cite{hausmann} stated (with a short outline of proof) 
that there is a bijective correspondence 
between $\Emb(\Sigma^3,S^6)$ and $C^3_3$. 
\end{remark}

Now recall that \textit{the Rohlin invariant} $\mu(\Sigma^3)$ of 
an integral homology $3$-sphere $\Sigma^3$ is defined, 
by choosing a smooth compact oriented spin $4$-manifold $W^4$ 
with $\partial{W^4}=\Sigma^3$, to be 
\[
\mu(\Sigma^3):=\sigma(W^4)/8\;\bmod2. 
\]
Then, we can characterise compressible embeddings $\Sigma^3\hookrightarrow S^6$ 
in the relation of our invariant $\Omega$ and the Rohlin invariant. 
Here, we say that an embedding $F\co\Sigma^3\hookrightarrow S^6$ is 
\textit{compressible} in $S^5$ if it is isotopic to 
an embedding in $S^5$ composed with the inclusion $j\co S^5\subset S^6$.  

\begin{corollary}\label{cor:compress1}
An embedding $F\co\Sigma^3\hookrightarrow S^6$ of 
an integral homology $3$-sphere $\Sigma^3$ in $S^6$ is 
compressible into $S^5$ if and only if 
$\Omega(F)\equiv\mu(\Sigma^3)\pmod2$. 
\end{corollary}

\begin{proof}
Assume that $F\co\Sigma^3\hookrightarrow S^6$ is isotopic to 
an embedding $f\co\Sigma^3\hookrightarrow S^5$. 
Let $\widetilde{f}\co V^4\hookrightarrow S^5$ be a Seifert surface for $f$. 
Then, $\widetilde{f}$ has trivial normal bundle and 
$V^4$ is necessarily a spin manifold (since it is embedded in codimension $1$). 
Therefore, we have $\Omega(F)=-(\sigma(V^4)-0)/8$, which is modulo $2$ equal to 
the Rohlin invariant $\mu(\Sigma^3)$. 

Conversely, assume that $\Omega(F)\equiv\mu(\Sigma^3)\pmod2$. 
By the same argument of Proposition~\ref{prop:conn}, 
$F$ is isotopic to $j\circ(\widetilde{f_0}|_{\Sigma^3}\conn g)$ 
for an embedding $\widetilde{f}_0\co W^4_0\hookrightarrow S^5$ 
of a spin $4$-manifold $W^4_0$ with $\partial{W^4_0}=\Sigma^3$ and 
an embedding 
$g\co S^3\hookrightarrow S^6$ of the $3$-sphere. 
Then, we have 
\begin{eqnarray*}
\Omega(F)
&=&\Omega(j\circ\widetilde{f_0}|_{\Sigma^3})+\Omega(g)\\
&\equiv&\mu(\Sigma^3)+\Omega(g)\bmod2.
\end{eqnarray*}
Thus, we have $\Omega(g)\equiv0\pmod2$, 
which implies 
by \cite[Theorem~5.17]{hae2} (see also \cite[Theorem~7.1]{takase2})
that $g\co S^3\hookrightarrow S^6$ 
can be isotoped to an embedding $g'$ of $S^3$ into $S^5$. 
Hence, $F$ is isotopic to 
$\widetilde{f}_0|_{\Sigma^3}\conn g'\co\Sigma^3\hookrightarrow S^6$ 
(composed with the inclusion $j$). 
\end{proof}

In view of the generalisation of Rohlin's theorem 
by Guillou and Marin \cite{g-m} (see also \cite{yukio2}), 
the following is a direct consequence of Corollary~\ref{cor:compress1}. 

\begin{corollary}\label{cor:compress2}
Let $\widetilde{F}\co V^4\hookrightarrow S^6$
be an embedding of a compact orientable $4$-manifold $V^4$ 
whose boundary $\partial{V^4}$ is diffeomorphic to 
an integral homology $3$-sphere and 
let $\Delta\subset V^4$ be an oriented surface 
dual to the normal Euler class of $\widetilde{F}$. 
Then, $\widetilde{F}|_{\partial{V^4}}\co\partial{V^4}\hookrightarrow S^6$ 
is compressible into $S^5$ if and only if 
$\mathrm{Arf}(V^4,\Delta)\equiv0\bmod2$.
\end{corollary}

\section{The homology bordism group}\label{sect:cobordism}
In this section, 
we show that the invariant $\Omega$ gives a homomorphism 
from the homology bordism group of embeddings of 
integral homology $3$-spheres in $S^6$. 

We say that two oriented integral homology $3$-spheres $\Sigma^3_0$ and $\Sigma^3_1$ 
are \textit{homology cobordant} if 
there exists a smooth compact oriented $4$-manifold $W^4$ 
with $\partial{W^4}=\Sigma^3_0\sqcup(-\Sigma^3_1)$ (the disjoint union) 
such that the inclusions induce isomorphisms 
$H_*(\Sigma^3_i)\to H_*(W^4)$ for $i=0,1$. 
Such a $W^4$ is called \textit{a homology cobordism} 
between $\Sigma^3_0$ and $\Sigma^3_1$. 
All homology cobordism classes of oriented integral 
homology $3$-spheres form an abelian group $\Theta^3_{\Z}$ 
with the group operation induced by the connected sum. 
The group $\Theta^3_{\Z}$ is infinite \cite{donaldson} and 
is infinitely generated \cite{furuta}. 

For embeddings of oriented integral 
homology $3$-spheres, we can naturally consider the notion of 
homology bordism (see \cite{saeki0}). 

\begin{definition}
Let $F_i$ $(i=0,1)$ be embeddings 
of oriented integral homology $3$-spheres $\Sigma^3_i$ in $S^6$ 
(or in $\R^6$). Then, $F_0$ and $F_1$ 
are said to be \textit{homology bordant} 
if there exist a homology cobordism $W^4$ 
between $\Sigma^3_0$ and $\Sigma^3_1$ and 
a proper embedding $H$ 
of $W^4$ in $S^6\times[0,1]$ (or in $\R^6\times[0,1]$) such that 
the embedding $H|_{\Sigma^3_i}$ of $\Sigma^3_i$ in $S^6\times\{i\}=S^6$ 
(or in $\R^6\times\{i\}=\R^6$) 
coincides with $F_i$ for $i=0,1$. 
Such an $H$ is called \textit{a homology bordism} between $F_0$ and $F_1$. 
\end{definition}

Now we prove that our invariant given in Proposition~\ref{prop:inv} 
(and in Theorem~\ref{thm:isotopic}) 
turns out to be invariant up to homology bordism. 

\begin{proposition}\label{prop:invcobor}
If two embeddings $F_i\co \Sigma^3_i\hookrightarrow S^6$ 
of \textup{(}distinct\textup{)} homology $3$-spheres $(i=0,1)$ 
are homology bordant, 
then $\Omega(F_0)=\Omega(F_1)$. 
\end{proposition}

\begin{proof}
We consider each embedding $F_i$ as an embedding 
$F_i\co \Sigma^3_i\hookrightarrow\R^6$ in $\R^6$ and let 
$H\co (W^4,\partial{W^4})\hookrightarrow (\R^6\times[0,1],\R^6\times\{0\}\sqcup\R^6\times\{1\})$
be a homology bordism between $F_0$ and $F_1$. 
We can assume that $H$ is perpendicular to $\R^6\times\{0\}\sqcup\R^6\times\{1\}$ 
along the boundary $\partial{W^4}$. 
Then, since $H^4(W^4)=H^3(W^4;\pi_2(\SO(3)))=H^2(W^4)=H^1(W^4)=0$, 
the normal bundle of $H$ is trivial and 
a normal framing $\nu=(\nu_1,\nu_3,\nu_3)$ for 
$H(W^4)\subset\R^6\times[0,1]$ 
determines normal framings for 
$F_i(\Sigma^3_i)\subset\R^6\times\{i\}$ ($i=0,1$). 

According to Lemma~\ref{lem:normal}, 
for each $i\in\{0,1\}$, 
consider a Seifert surface 
$\widetilde{F}_i\co V^4_i\hookrightarrow\R^6$ 
for $F_i$ such that the outward normal vector field of 
$F(\Sigma^3_i)\subset\widetilde{F}_i(V^4_i)$ coincides with 
the first vector field $\nu_1$ of $\nu$. 

\begin{figure}[tbp]
\begin{center}
\includegraphics[height=30mm,keepaspectratio,bb=196 569 341 718]{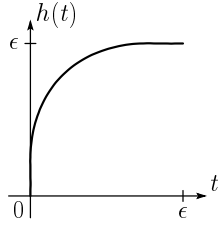}
\end{center}
\caption{``Smoothing function $h$''}\label{fig:func}
\end{figure}

Take collar neighbourhoods 
$C_i:=\Sigma^3_i\times[0,\epsilon]$ of 
the boundary $\partial{V^4_i}=\Sigma^3_i$. 
Then, by using $\widetilde{F}_0$, $\widetilde{F}_1$, $H$ and 
an appropriate ``smoothing function'' 
$h\co[0,\epsilon]\to[0,\epsilon]$ 
as in Figure~\ref{fig:func}, 
we can define a smooth embedding 
\[
G\co V^4_0\cup W^4\cup (-V^4_1)\hookrightarrow\R^6\times[-\epsilon,1+\epsilon]\subset\R^7
\]
of the closed $4$-manifold $M^4:=V^4_0\cup W^4\cup (-V^4_1)$ by 
\[G(x)=\left\{\begin{array}{lcl}
(\widetilde{F}_0(x),-\epsilon)&\text{ if }&x\in V^4_0\smallsetminus C_0\\
(\widetilde{F}_0(\bar{x},t),-h(t))&\text{ if }&(\bar{x},t)\in\Sigma^3_0\times[0,\epsilon]=C_0\\
H(x)&\text{ if }&x\in W^4\\
(\widetilde{F}_1(\bar{x},t),1+h(t))&\text{ if }&(\bar{x},t)\in\Sigma^3_1\times[0,\epsilon]=C_1\\
(\widetilde{F}_1(x),1+\epsilon)&\text{ if }&x\in V^4_1\smallsetminus C_1. 
\end{array}\right.\]

Furthermore, 
with $\theta_t\in[0,\pi/2]$ such that $h'(t)=\tan\theta_t$, 
we can define a non-zero normal vector field $n$ of $G(M^4)\subset\R^7=\{(y_1,y_2,\cdots,y_7)\}$ as: 
\[\left\{\begin{array}{lcl}
-{\partial}/{\partial{y_7}}&\text{ on }&G(V^4_0\smallsetminus C_0)\\
-\cos\theta_t\cdot({\partial}/{\partial{y_7}})+\sin\theta_t\cdot\nu_1&\text{ on }&G(*,t)\text{ where }(*,t)\in\Sigma^3_0\times[0,\epsilon]=C_0\\
\nu_1&\text{ on }&W^4\\
\cos\theta_t\cdot({\partial}/{\partial{y_7}})+\sin\theta_t\cdot\nu_1&\text{ on }&G(*,t)\text{ where }(*,t)\in\Sigma^3_1\times[0,\epsilon]=C_1\\
{\partial}/{\partial{y_7}}&\text{ on }&G(V^4_1\smallsetminus C_1). 
\end{array}\right.\]

Denote the push-off of $G(M^4)$ 
along the normal field $n$ by $G(M^4)'$. 
Then, for any $c=(c_1,c_2)\in H_2(M^4)\approx H_2(V_0)\oplus H_2(-V_1)$, 
the linking number $\lk(c,G(M^4)')$ in $\R^7$ 
is equal to $\lk(c_1,G(M^4)')+\lk(c_2,G(M^4)')$, 
which vanishes since $c_0$ and $c_1$ bound $3$-chains respectively 
in $\R^6\times\{-\epsilon\}$ and in $\R^6\times\{1+\epsilon\}$. 
Therefore, 
the normal vector field $n$ is unlinked (\textit{non enlac\'e}) 
in the sense of Bo\'echat and Haefliger \cite[Definition~1.1]{b-h}. 

Now, we want to apply Bo\'echat and Haefliger's argument \cite{b-h} to 
the embedding $G$ with the normal field $n$. 
Let $e$ be the Euler class of the $2$-plane bundle 
complementary to $n$ in the normal bundle of $G(M^4)\subset\R^7$. 
Then, by \cite[Th\'eor\`eme 2.1]{b-h}, we have 
\[
\sigma(M^4)-e\smile e=0. 
\]
Since $\sigma(M^4)=\sigma(V_0)-\sigma(V_1)$ and 
$e=e_{\widetilde{F}_0}-e_{\widetilde{F}_1}\in H_2(M^4)\approx H_2(V_0)\oplus H_2(-V_1)$, 
\[
0=(\sigma(V_0)-\sigma(V_1))-(e_{\widetilde{F}_0}\smile e_{\widetilde{F}_0}-e_{\widetilde{F}_1}\smile e_{\widetilde{F}_1}). 
\]
Thus we have 
$\sigma(V_0)-e_{\widetilde{F}_0}\smile e_{\widetilde{F}_0}=\sigma(V_1)-e_{\widetilde{F}_1}\smile e_{\widetilde{F}_1}$. 
This completes the proof. 
\end{proof}

\begin{remark}
In the above proof of Proposition~\ref{prop:invcobor}, 
the Seifert surfaces $\widetilde{F}_i\co V^4_i\hookrightarrow\R^6$ ($i=0,1$) 
are chosen in such a way that the normal field $n$ is 
actually non-zero on the whole $G(M^4)$, 
so that we actually have also $3\sigma(M^4)+e\smile e=0$, 
and hence $\sigma(M^4)=e\smile e=0$ (see \cite[p94]{g-m}). 
\end{remark}

All homology bordism classes of embeddings 
of oriented integral homology $3$-spheres in $S^6$ 
form an abelian group, denoted by $\Gamma^3_3$, 
with group structure via connected sum. 
Then, Proposition~\ref{prop:invcobor} implies that 
$\Omega$ induces the homomorphism 
$\Omega\co\Gamma^3_3\to\Z$, 
since Hopf invariants and signatures of Seifert surfaces 
are additive under boundary connected sum 
(by Proposition~\ref{prop:additivity} and by Novikov additivity). 
Furthermore, the invariant $\Omega$ leads to 
the following characterisation of $\Gamma^3_3$. 

\begin{theorem}\label{thm:cobordism}
The following is an isomorphism: 
\[\begin{array}{cccc}
\bar{\Omega}\co&\Gamma^3_3&\stackrel{\approx}{\to}&\Theta^3_{\Z}\,\oplus\,\Z\\
&{[F\co\Sigma^3\hookrightarrow S^6]}&\mapsto&(\,[\Sigma^3],\,\Omega(F)\,), 
\end{array}\]
where $[\Sigma^3]$ is the homology cobordism class represented 
by $\Sigma^3$ and $[F\co\Sigma^3\hookrightarrow S^6]$ is  
the homology bordism class represented 
by $F\co\Sigma^3\hookrightarrow S^6$. 
\end{theorem}

\begin{proof}
Since any integral homology $3$-sphere can be embedded in $S^6$, 
the surjectivity follows from Theorem~\ref{thm:3dim} and 
the additivity of $\Omega$ under connected sum. 

Now, let $\Sigma^3$ be an integral homology sphere 
homology cobordant to zero (i.e., bounding an bounding acyclic 4-manifold $V^4$) 
and $F\co\Sigma^3\hookrightarrow S^6$ be 
an embedding with $\Omega(F)=0$. 
Then, $V^4$ can be embedded in $S^5$; 
let $f\co V^4\hookrightarrow S^5$ be an embedding. 
Clearly, 
$\Omega((j\circ f)|_{\partial{V^4}})=-\sigma(V^4)/8=0$, 
where $j\co S^5\subset S^6$ is the inclusion.  
Therefore, by Theorem~\ref{thm:isotopic}, 
$F$ is isotopic to $(j\circ f)|_{\partial{V^4}}$ 
and hence is homology bordant to zero. 
Since $\bar{\Omega}$ is clearly a homomorphism, 
it should be an isomorphism. 
\end{proof}

\begin{corollary}\label{cor:cobordism}
Two embeddings of an oriented integral homology 
$3$-sphere in $S^6$ are isotopic if and only if 
they are homology bordant. 
\end{corollary}

\begin{acknowledgement}
The author would like to thank 
Professor Dennis Roseman for many helpful comments. 
The author is partially supported by the Grant-in-Aid for JSPS Fellows.
\end{acknowledgement}


\begin{thebibliography}{00}


\bibitem{boechat}J.~Bo\'echat, 
{\em Plongements de vari\'et\'es diff\'erentiables orient\'ees de dimension $4k$ dans $\mathbf{R}^{6k+1}$},
Comment. Math. Helv. \textrm{46} (1971) 141--161. 


\bibitem{b-h}J.~Bo\'echat and A.~Haefliger, 
{\em Plongements diff\'erentiables des vari\'et\'es orient\'ees de dimension $4$ dans $\mathbf{R}^7$},
Essays on Topology and Related Topics, 
M\'emoires d\'edi\'es \`a Georges de Rham, 156--166, Springer, New York, 1970. 


\bibitem{donaldson}S.~K.~Donaldson, 
{\em An application of gauge theory to four-dimensional topology},
J.\ Differential Geom.\ \textrm{18} (1983) 279--315. 


%
%


%
%
\bibitem{freedman}M.~Freedman, 
{\em The topology of four-dimensional manifolds},
J.\ Differential Geom.\ \textrm{17} (1982) 357--453. 


\bibitem{furuta}M.~Furuta, 
{\em  Homology cobordism group of homology $3$-spheres},
Invent.\ Math.\ \textrm{100} (1990) 339--355.


\bibitem{g-m}L.~Guillou and A.~Marin, 
{\em Commentaires sur les quatre articles pr\'ec\'edents de V. A. Rohlin},
in: L.~Guillou, A.~Marin (Eds.),
\`A la recherche de la topologie perdue, 25--95, 
Progr. Math., 62, Birkh\"auser Boston, Boston, MA, 1986. 


\bibitem{gompf}R.~E.~Gompf, 
{\em Handlebody construction of Stein surfaces},
Ann.\ of Math. \textrm{148} (1998) 619--693.


\bibitem{hae1}A.~Haefliger, 
{\em Knotted $(4k-1)$-spheres in $6k$-space},
Ann.\ of Math. \textrm{75} (1962) 452--466.


\bibitem{hae2}A.~Haefliger, 
{\em Differential embeddings of $S^n$ in $S^{n+q}$ for $q\ge2$},
Ann.\ of Math. \textrm{83} (1966) 402--436.


\bibitem{hausmann}J.~Hausmann, 
{\em Plongements de sph\`eres d'homologie},
C.\ R.\ Acad.\ Sci.\ Paris S\'er.\ A-B \textrm{275} (1972), A963--965.


%
%
\bibitem{hirsch2}M.~Hirsch, 
{\em Embeddings and compressions of polyhedra and smooth manifolds},
Topology \textrm{4} (1966) 361--369.


%
%


\bibitem{kirby}R.~Kirby, 
{\em The topology of $4$-manifolds}, 
Lecture Notes in Mathematics \textrm{1374}, 
Springer-Verlag, Berlin, 1989.


\bibitem{kuperberg}G.~Kuperberg, 
{\em Noninvolutory Hopf algebras and $3$-manifold invariants},
Duke Math.\ J.\ \textrm{84} (1996) 83--129.


\bibitem{levine}J.~Levine, 
{\em Imbedding and isotopy of spheres in manifolds}, 
Proc.\ Camb.\ Phil.\ Soc.\ \textrm{60} (1964) 433--437. 


\bibitem{levine1}J.~Levine, 
{\em Semi-free circle actions on spheres}, 
Invent.\ Math.\ \textrm{22} (1973) 161--186. 


\bibitem{lu}Z.~L\"u, 
{\em  The generalized Smith conjecture of codimension greater than two},
J.\ Knot Theory Ramifications \textrm{9} (2000), 479--490.


%
%
%


\bibitem{masuda}M.~Masuda, 
{\em  Smooth involutions on homotopy $C\mathrm{P}^3$},
Amer.\ J.\ of Math. \textrm{106} (1984) 1487--1501.




\bibitem{yukio2}Y.~Matsumoto, 
{\em  An elementary proof of Rochlin's signature theorem and its extension by Guillou and Marin},
\`A la recherche de la topologie perdue, 119--139, 
Progr. Math., 62, Birkh\"auser Boston, Boston, MA, 1986. 


\bibitem{m-y}D.~Montgomery and C.~T.~Yang, 
{\em  Differentiable actions on homotopy seven spheres},
Trans.\ Amer.\ Math.\ Soc.\ \textrm{122} (1966) 480--498. 
%
%


\bibitem{pontrjagin}L.~Pontrjagin, 
{\em A classification of mappings of the three-dimensional complex into the two-dimensional sphere},
Rec.\ Math.\ [Mat. Sbornik] N.\ S.\ \textrm{9 (51)} (1941) 331--363.


\bibitem{saeki0}O.~Saeki, 
{\em Cobordism classification of knotted homology $3$-spheres in $S^5$},
Osaka J.\ Math.\ \textrm{25} (1988) 213--222.


\bibitem{saeki}O.~Saeki, 
{\em On punctured $3$-manifolds in $5$-sphere},
Hiroshima Math.\ J.\ \textrm{29} (1999) 255--272.

%
%
%
%


\bibitem{smale2}S.~Smale, 
{\em On the structure of manifolds},
Amer.\ J.\ of Math. \textrm{84} (1962) 387--399.


%
%


%
%
\bibitem{takase2}M.~Takase, 
{\em A geometric formula for Haefliger knots},
Topology \textrm{43} (2004) 1425--1447.


\bibitem{takase3}M.~Takase, 
{\em The Hopf invariant of a Haefliger knot},
Math.\ Z.\ (to appear), (arXiv: math.GT/0506356).


\bibitem{u-m}H.~Uehara and W.~S.~Massey, 
{\em The Jacobi identity for Whitehead Products}, 
Algebraic geometry and topology. 
A symposium in honor of S. Lefschetz, 361--377, 
Princeton University Press, Princeton, N. J., 1957. 


\bibitem{white}J.~White, 
{\normalfont Twist invariants and the Pontrjagin numbers of immersed manifolds},
Differential geometry (Proc. Sympos. Pure Math., Vol. XXVII, Part 1, 
Stanford Univ., Stanford, Calif., 1973), 429--437, 
Amer. Math. Soc., Providence, R. I., 1975. 




\end{thebibliography}
\end{document}